\documentclass[draft]{amsart}
\usepackage{amsmath,amssymb,amsthm,latexsym}
\usepackage{url}
\usepackage{amscd}
\theoremstyle{definition}
\newtheorem{thm}{Theorem}[section]
\newtheorem{lem}[thm]{Lemma}
\newtheorem{cor}[thm]{Corollary}
\newtheorem{prop}[thm]{Proposition}
\newtheorem{defn}{Definition}[section]

\newtheorem{rem}{Remark}[section]

\newtheorem{ques}{Question}[section]

\DeclareMathOperator{\tr}{Tr}
\DeclareMathOperator{\diag}{diag}
\DeclareMathOperator*{\Max}{Max}
\DeclareMathOperator{\rank}{rank}
\newcommand{\mt}{\otimes_{\min}}
\begin{document}
\date{}
\title{Hilbertian matrix cross normed spaces arising from normed ideals}
\author{Takahiro Ohta}

\address{Department of Mathematics, Graduate School of Science, Kyoto University, Kyoto 606-8502, Japan}
\email{ohta@math.kyoto-u.ac.jp}

\begin{abstract}
Generalizing Pisier's idea, we introduce a Hilbertian matrix cross normed space associated with
a pair of symmetric normed ideals. When the two ideals coincide, we show that our construction
gives an operator space if and only if the ideal is the Schatten class.
In general, a pair of symmetric normed ideals that are not necessarily the Schatten class 
may give rise to an operator space.
We study the space of completely bounded mappings between the matrix cross normed spaces
obtained in this way and show that the multiplicator norm naturally appears as the completely bounded norm.
\end{abstract}

\maketitle

\section{Introduction}\label{intro}

An operator space is a subspace of the set of bounded operators on a Hilbert space,
which is abstractly characterized as a Banach space equipped with matrix norms satisfying certain properties.
An operator space whose base space is a Hilbert space 
is said to be a Hilbertian operator space. 
The theory of homogeneous Hilbertian operator space is one of the 
central topics in operator space theory and it plays an essential role in various situations.
For example, it is used to analyze the structures of the space of operator spaces
with the metric which is analogous to the Banach-Mazur distance (cf.~\cite{P2})
and to obtain an embedding of operator spaces into noncommutative $L_{p}$-spaces (cf.~\cite{J} and \cite{P3}).

The relationships between homogeneous Hilbertian operator spaces and 
operator ideals are first studied by Mathes and Paulsen.
Mathes and Paulsen considered in \cite{MP} a larger category, called 
matricially normed spaces (m.c.n.\ spaces), than that of operator spaces. 
They showed that if $H_{1}$ and $H_{2}$ are homogeneous Hilbertian m.c.n.\ spaces with 
the common base space $H$, then the space of completely bounded mappings
$CB(H_{1},H_{2})$ becomes a symmetric normed 
ideal (s.n.\ ideal) \cite[1.2.\ Proposition]{MP} and showed 
that every s.n.\ ideal on $B(H)$ which is not 
equivalent to the ideal of compact operators or the ideal of trace class operators is 
isomorphic as a set to the space of completely bounded mappings 
on some homogeneous Hilbertian m.c.n.\ spaces \cite[2.2.\ Theorem]{MP}.

G. Pisier showed that the norm of the elements in the interpolating spaces 
between the row Hilbert space and the column Hilbert space
is represented by the operator norm on the Schatten ideals \cite[Theorem 8.4]{P2}.
Inspired by this analysis, 
in our paper we introduce a Hilbertian m.c.n.\ space $H(\Phi,\Psi)$ for a 
pair of symmetric norming functions (s.n.\ functions) $\Phi, \Psi$ with 
$\Phi\geq\Psi$ and investigate the structure of the space. 
The matrix norm of $H(\Phi,\Psi)$ is defined by
\[ \| T\|_{H(\Phi,\Psi)}
=\left(\sup_{x}\frac{\|\sum T_{i}x T_{i}^{\ast}\|_{\Psi}}{\| x\|_{\Phi}}\right)^{1/2}, \]
where $T=\sum\xi_{i}\otimes T_{i}\in H\otimes M_{n}$ and 
$(\xi_{i})$ is an orthonormal basis of a separable Hilbert space $H$.
We also focus on the space of 
completely bounded mappings between two spaces arising in this way.
The m.c.n.\ space $H(\Phi,\Psi)$ is not always an operator space.
In section \ref{mcn} we show that if the m.c.n.\ space $H(\Phi,\Psi)$ is an operator space,
then for all $x, y, z\in\mathfrak{S}_{\Phi}$ the following inequality
\[ \frac{\| x\otimes y\|_{\Psi}}{\| x\|_{\Psi}}
\leq\frac{\| z\otimes y\|_{\Phi}}{\| z\|_{\Phi}} \]
is satisfied, where $\mathfrak{S}_{\Phi}$ is the s.n.\ ideal
arising from $\Phi$.
In particular, when $\Phi=\Psi$ we show that
the m.c.n.\ space $H(\Phi)=H(\Phi,\Phi)$ is an operator space 
if and only if $\Phi$ is the Schatten norm. 
However, the situation differs for $\Phi\neq\Psi$.
Indeed, when $\Phi$ is a Q$^{\ast}$-norm 
and $\Psi$ is a Q-norm, 
$H(\Phi,\Psi)$ is always an operator space.

We also study the space of completely bounded mappings 
between m.c.n.\ spaces we constructed.
We determine the completely bounded norm from 
the row Hilbert space $R$ to $H(\Phi,\Psi)$ as
\[ \| x\|_{CB(R, H(\Phi,\Psi))}=
\left(\sup_y\frac{\big\||x|^{2}\otimes y\big\|_{\Psi}}
{\| y\|_{\Phi}}\right)^{1/2}. \]
This implies that if $H(\Phi,\Psi)$ is an operator space, 
then we have the isometric isomorphisms
$CB(R, H(\Phi,\Psi))=\mathfrak{S}_{\tilde{\Psi}}$ and 
$CB(C, H(\Phi,\Psi))=\mathfrak{S}_{\tilde{\Phi^{\ast}}}$
for the column Hilbert space $C$ 
(see section \ref{mcn} for the definition of $\tilde{\Phi}$).

The above result leads us to consider the condition:
\[ \exists c>0,~
 \| x\otimes y\|_{\Psi}\leq c\| x\|_{\Psi}\| y\|_{\Psi},~\forall x,y\in\mathfrak{S}_{\Psi}. \]
This condition implies that there exists a constant
\[p=\lim_{n\rightarrow\infty}\frac{\log n}
{\log\| P_{n}\|_{\Phi}}\quad (P_{n} \text{ is any rank } n \text{ projection})\]
such that $\| x\|_{p}\leq c\| x\|_{\Phi}$,
where $\| x\|_{p}$ is the Schatten $p$-norm.
This together with a dual version implies
the above mentioned fact that $H(\Phi)$ is an operator space only if
$\Phi$ is the Schatten norm.


\section{Preliminaries}
In this section we collect the basics of the theory of operator spaces and operator ideals,
which are often used in the paper.
We refer to \cite{ER} and \cite{P1} for the theory of operator spaces and
to \cite{GK} for the theory of operator ideals.

An operator space is abstractly characterized as follows.
We consider a Banach space $E$ such that for 
each $n\in\mathbb{N}$ there is a norm $\|\cdot\|_{n}$ on the 
matrix space $M_{n}(E)$ of $n\times n$ matrices with entries in the elements of 
$E$ and the family $\{M_{n}(E),~\|\cdot\|_{n}\}$ with $\|\cdot\|_{1}$ equal to the 
original norm of $E$.
Then we can consider the two properties
\begin{enumerate}
\item[(M1)] $\left\|\begin{pmatrix}
    x & 0\\
    0 & y
    \end{pmatrix}\right\|_{m+n}=\Max\{\| x\|_{m},\| y\|_{n}\}$ for any 
$x\in M_{m}(E),~y\in M_{n}(E)$, and $m,n\in\mathbb{N}$, and
\item[(M2)] $\| a x b\|_{n}\leq\| a\|\| x\|_{m}\| b\|$ for any $x\in M_{m}(E)$,~ 
$a\in M_{n\times m},~b\in M_{m\times n}$, and $m,n\in\mathbb{N}$, where 
$M_{m\times n}=M_{m\times n}(\mathbb{C})$ and $a x b$ 
means the matrix product. 
\end{enumerate}
(M1) may be replaced with
\begin{enumerate}
\item[(M1)$^{\prime}$] $\left\|\begin{pmatrix}
    x & 0\\
    0 & y
    \end{pmatrix}\right\|_{m+n}\leq\Max\{\| x\|_{m},\| y\|_{n}\}$, for any 
$x\in M_{m}(E),y\in M_{n}(E)$, and $m,n\in\mathbb{N}$.
\end{enumerate}
For a Hilbert space $H$ an operator space $E\subseteq B(H)$ is a Banach space 
satisfying the 
properties (M1) and (M2) under the identification of $M_{n}(E)$ as a 
subspace of $M_{n}(B(H))=B(H^{n})$. 
Conversely, Ruan~\cite[Theorem 3.1]{R} showed that 
a Banach space having the matrix norm structure with the 
properties (M1) and (M2) has an isometric embedding into 
the space $B(H)$ for some Hilbert space $H$ such that the matrix norms
come from $M_{n}(B(H))=B(H^{n})$.
The properties (M1) and (M2) are called Ruan's axioms. 
In the operator space category, the morphisms are 
the completely bounded (c.b.) mappings. Let $E$, $F$ be operator spaces and $u$ be a 
linear mapping from $E$ to $F$. We say that $u$ is completely bounded if
\[\| u\|_{cb}=\sup_{n}\| id_{n}\otimes u\colon M_{n}(E)\rightarrow M_{n}(F)\| <\infty,\]
where $M_{n}(E)$ is identified with the algebraic tensor product $M_{n}\otimes E$.
The completely bounded norm of $u$ is defined by $\| u\|_{cb}$.
An operator space $E$ is said to be
homogeneous if for any bounded linear mapping $u$ on $E$ we have 
$\| u\|=\| u\|_{cb}$.
We denote the Banach space of completely bounded mappings from $E$ to $F$ 
with norm $\|\cdot\|_{cb}$ by $CB(E,F)$.

The category of matrix cross normed spaces is larger than 
that of operator spaces.
Let $H$ be a separable Hilbert space with a sequence of matrix norms $\{\|\cdot\|_{n}\}_{n=1}^{\infty}$
on the family $\{ M_{n}(H)\}_{n=1}^{\infty}$ such that $\|\cdot\|_{1}$ coincides with the norm of $H$.
We call $H$ a matrix cross normed space (m.c.n.\ space) if
\[\| x\otimes A\|_{n}=\| x\|\| A\|_{M_{n}} \]
for all $x\in H,~A\in M_{n}$, and $n\in\mathbb{N}$.

For a finite-dimensional or separable infinite-dimensional Hilbert space 
$K$ with dimension $n$, identifying $B(K)$ with the matrix space $M_{n}$ 
we denote the matrix whose 
$(i,j)$-entry is 1 and the other entries are $0$ by $e_{i j}$.

Next we introduce the basic theory of the operator ideals (cf. \cite[Chapter III]{GK}).
Let $c_0$, $\hat{c}$, and $\hat{k}$ be the spaces of sequences of real numbers 
defined by
\begin{enumerate}
\item[]$c_0=\left\{\xi=\{\xi_{i}\}:\lim_{i\rightarrow\infty}\xi_{i}=0\right\},$
\item[]$\hat{c}=\big\{\xi=\{\xi_{i}\}\in c_0:
    \text{only finitely many}~\xi_{i}\text{'s are nonzero}\big\},$
\item[]$\hat{k}=\big\{\xi=\{\xi_{i}\}\in\hat{c}:
    \xi_{1}\geq\xi_{2}\geq\ldots\xi_{n}\geq\ldots\geq 0\big\},$
\end{enumerate}
respectively. A real valued function $\Phi$ on $\hat{c}$ is 
called a symmetric norming (s.n.) function if it satisfies the followings:
\begin{enumerate}
\item[(1)]$\Phi$ is a norm on $\hat{c}$;
\item[(2)]$\Phi(1,0,0,\ldots)=1$;
\item[(3)]$\Phi(\xi_{1},\xi_{2},\ldots,\xi_{n},0,0,\ldots)=
 \Phi(|\xi_{j_{1}}|,|\xi_{j_{2}}|,\ldots,|\xi_{j_{n}}|,0,0,\ldots)$
 for all $\xi\in\hat{c}$,
where $\{j_{1},j_{2},\ldots,j_{n}\}$ is any permutation of $\{1,2,\ldots,n\}$.
\end{enumerate}
For an s.n.\ function $\Phi$, we set
\[c_{\Phi}=\big\{\xi=\{\xi_{i}\}\in c_0:\sup_{n}\Phi(\xi^{(n)})<\infty\big\},\]
where $\xi^{(n)}=(\xi_{1},\ldots,\xi_{n},0,0,\ldots)$.
We extend the domain of $\Phi$ by
\[\Phi(\xi)=\lim_{n\rightarrow\infty}\Phi(\xi^{(n)}),~\xi\in c_{\Phi}.\]
For $1\leq p \leq\infty$, we denote by $\Phi_{p}$ the $\ell_{p}$-norm.

Throughout the paper, $H$ denotes a separable infinite-dimensional Hilbert 
space with an orthonormal basis $\{\xi_{i}\}_{i=1}^{\infty}$ and 
$\mathfrak{S}_{\infty}$ denotes the subspace of $B(H)$ 
consisting of all compact operators on $H$.
For $x\in\mathfrak{S}_{\infty}$ we denote by $\{s_{j}(x)\}_{j=1}^{\infty}$ the 
singular numbers ($s$-numbers) of $x$,
i.e.\ the nonincreasing rearrangement of eigenvalues of $|x|$.

Let $\mathfrak{S}$ be a two-sided ideal of $B(H)$. A functional 
$\|\cdot\|_s$ on $\mathfrak{S}$ is said to be a symmetric norm if 
it satisfies the followings:
\begin{enumerate}
\item[(1)]
$\|\cdot\|_s\text{ is a norm on }\mathfrak{S};$
\item[(2)]
$\text{for any rank one operator }x,~\| x\|_s=\| x\|;$
\item[(3)]
$\| a x b\|_s\leq\| a\|\| x\|_s\| b\|~(\forall a,b\in B(H),~
    \forall x\in\mathfrak{S}).$
\end{enumerate}
We call $(\mathfrak{S},\|\cdot\|_s)$ a symmetrically normed ideal if 
$\|\cdot\|_s$ is a symmetric norm on $\mathfrak{S}$
and makes $\mathfrak{S}$ a Banach space.

For an s.n.\ function $\Phi$, 
we denote by $\mathfrak{S}_{\Phi}$ the set of operators 
$x\in\mathfrak{S}_{\infty}$ with $s(x)=\{s_{j}(x)\}\in c_{\Phi}$, and put
\[\|x\|_{\Phi}=\Phi(s(x)).\]
Then $\mathfrak{S}_{\Phi}$ is an s.n.\ ideal with the norm $\|\cdot\|_{\Phi}$.
In this paper we often use the property
\[xx^{\ast}\in\mathfrak{S}_{\Phi}\Leftrightarrow 
x^{\ast}x\in\mathfrak{S}_{\Phi} \text{ and } 
\| xx^{\ast}\|_{\Phi}=\| x^{\ast}x\|_{\Phi}.\]

Let $\Phi$ be an s.n.\ function. The function 
\[\Phi^{\ast}(\eta)=\max_{\xi\in\hat{k}}\left\{\frac{1}{\Phi(\xi)}
 \sum_{i}\eta_{i}^{\ast}\xi_{i}\right\}.\]
makes sense for any $\eta\in\hat{c}$ and $\Phi^{\ast}$ is an s.n.\ function. 
We call $\Phi^{\ast}$ the 
adjoint of $\Phi$. Note that for any s.n.\ function $\Phi$, we have  
$(\Phi^{\ast})^{\ast}=\Phi$ and the following duality 
\[\| x\|_{\Phi}=\sup_{\| y\|_{\Phi^{\ast}}\leq 1}|\tr(y x)|.\]

We introduce a few classes of normed ideals used in this paper.
We denote by $\mathfrak{S}_{p}=\mathfrak{S}_{\Phi_{p}}$
the Schatten ideal for $1\leq p\leq\infty$.
For $1\leq q\leq p<\infty$, the Lorentz ideal $S_{p,q}$
is an s.n.\ ideal whose norm is given by
\[\| x\|_{p,q}=\left(\sum_{j=1}^{\infty}\frac{s_{j}(x)^q}{j^{1-q/p}}\right)^
{1/q}.\]
Let $1=\pi_{1}\geq\pi_{2}\geq\cdots\geq0$ be
a sequence of nonincreasing positive numbers such that
$\lim_{n\rightarrow\infty}\pi_{n}=0$ and
$\sum_{n=1}^{\infty}\pi_{n}=\infty$.
We say that such a sequence is binormalizing.
The s.n.\ function $\Phi_{\pi}$ is defined by 
\[\Phi_{\pi}(a)=\sum_{n=1}^{\infty}\pi_{n}a_{n}^{\ast},~a=(a_{n}),\]
where $(a_{n}^{\ast})$ is the nonincreasing rearrangement of $(a_{n})$.
Note that if $q=1$, then the Lorentz ideal $S_{p,1}$ is equal to the ideal $\mathfrak{S}_{\Phi_{\pi}}$ 
defined by the binormalizing sequence $\pi_{j}=j^{1/p-1}$.

Finally we introduce an important class of operator spaces.
If $E_{0},E_{1}$ are compatible Banach spaces, then we denote by $(E_{0},E_{1})_{\theta}$ 
for $0<\theta <1$ the complex interpolation space of them (see~\cite[Chapter 4]{BL}). 
If $E_{0},E_{1}$ are operator spaces whose base spaces are compatible, we construct 
an operator space complex interpolation by identifying $M_{n}((E_{0},E_{1})_{\theta})$ 
with $(M_{n}(E_{0}),M_{n}(E_{1}))_{\theta}$ for each $n\in\mathbb{N}$.
We denote by $R$ and $C$ the row and column operator space respectively \cite[Section 3.4]{ER}.
These spaces are homogeneous Hilbertian operator spaces whose matrix norms are given by
\[\left\| \sum_{i=1}^{n}\xi_{i}\otimes T_{i}\right\|_{R}=
\left\|\sum_{i=1}^{n} T_{i}T_{i}^{\ast}\right\|^{1/2},~
\left\| \sum_{i=1}^{n}\xi_{i}\otimes T_{i}\right\|_{C}=
 \left\|\sum_{i=1}^{n} T_{i}^{\ast}T_{i}\right\|^{1/2},\]
for a finite sequence of matrices $\{T_{i}\}_{i=1}^{n}$.
Note that $R^{\ast}=C$ and $C^{\ast}=R$ in the operator space category.
We denote by $R(\theta)$ the operator space complex 
interpolation $(R,C)_{\theta}$ for $0<\theta <1$,
which is a homogeneous Hilbertian operator space. 
We set $R(0)$ to be the row Hilbert space $R$ and 
$R(1)$ to be the column Hilbert space $C$.
When $\theta=1/2$, we write $OH=R(1/2)$. 
Pisier~\cite[Theorem 1.1]{P2} introduced these spaces and showed that 
for any finite sequence $\{T_{i}\}$ it holds that
\[\left\| \sum_{i}\xi_{i}\otimes T_{i}\right\|_{OH}=
\left\|\sum_{i} T_{i}\otimes\bar{T_{i}}\right\|^{1/2},\]
where $\bar{T_{i}}$ means the complex conjugate of $T_{i}$.
Another important property of $OH$ is the self-duality. 
For an operator space $E$, the operator space $\bar{E}$ means its complex conjugate. 
The matrix norms of the elements of $\bar{E}$ are defined by 
\[ \| (\overline{x_{i j}})\|_{M_{n}(\bar{E})}=\| (x_{i j})\|_{M_{n}(E)}. \]
Pisier showed in \cite[Theorem 1.1]{P2} the completely isometric identification
\[ OH=\overline{OH^{\ast}}. \]
Another important example of a homogeneous Hilbertian operator space is the 
minimal operator space $H_{\min}$.
Let $E$ be a Banach space. 
We can embed $E$ into a commutative $C^{\ast}$-algebra 
(for example the space of all continuous functions 
on the unit ball of $E^{\ast}$ equipped with the weak topology). 
We denote by $\min(E)$ the operator space whose matrix norms arise form this embedding. 
The minimal operator space norm is the minimal norm among all operator space norms. 
When $E$ is a Hilbert space $H$, we denote the minimal operator space by $H_{\min}$.
The matrix norm on $H_{\min}$ satisfies 
\[\left\| \sum_{i=1}^{m}\xi_{i}\otimes T_{i}\right\|_{\min}=
\sup\left\|\sum_{i=1}^{m}v_{i}T_{i}\right\|,\]
where the supremum is taken over all unit vectors $\{v_{i}\}$ of $\ell_{2}^{m}$.


\section{Basic Properties of the m.c.n.\ space $H(\Phi,\Psi)$}\label{mcn}

Let $K$ be a separable Hilbert space which is identified with a subspace of separable 
infinite-dimensional Hilbert space. 
For $n\in\mathbb{N}\cup\{\infty\}$ we denote by $I_{n}$ the identity operator 
on the Hilbert space of dimension $n$.
Let $T$ be a finite sum $T=\sum_{i}\xi_{i}\otimes T_{i}$
in the algebraic tensor product $H\otimes B(K)$ and
we set $T^{\ast}=\sum_{i}\xi_{i}\otimes T_{i}^{\ast}$.
Pisier showed the identification of matrix norms of $R(\theta)~(0\leq \theta\leq 1)$ in \cite[Theorem 8.4]{P2} as follows:
\[ \left\|\sum_{i}\xi_{i}\otimes T_{i}\right\|_{R(\theta)\mt B(K)}
=\sup\left\{\left\|\sum_{i}T_{i}x T_{i}^{\ast}\right\|_{p}^{1/2}:x\in\mathfrak{S}_{p,+},~
\| x\|_{p}\leq 1\right\}, \]
where $p=\theta^{-1}$.
We define the operators 
$\rho_{T}$ and $\rho_{T^{\ast}}$ on $B(K)$ by
\[\rho_{T}(x)=\sum T_{i}x T_{i}^{\ast},~x\in B(K),\]
\[\rho_{T^{\ast}}(x)=\sum T_{i}^{\ast}x T_{i},~x\in B(K).\]
Neither $\rho_{T}$ nor $\rho_{T^{\ast}}$ depends on the choice 
of the basis $\{\xi_{i}\}_{i=1}^{\infty}.$
If $\mathfrak{S}$ is a two-sided ideal in $B(K)$, we have $\rho_{T}(\mathfrak{S})
\subseteq \mathfrak{S}$ and $\rho_{T^{\ast}}(\mathfrak{S})\subseteq\mathfrak{S}$.
For fixed s.n.\ functions $\Phi$ and $\Psi$ with $\Psi\leq\Phi$, we define a norm
$\|\cdot\|_{\Phi,\Psi}$ 
on the space of finite sums $T\in H\otimes B(K)$ by
\[ \| T\|_{\Phi,\Psi}=\|\rho_{T}\colon\mathfrak{S}_{\Phi}\rightarrow\mathfrak{S}_{\Psi}\|^
{1/2}. \]
Now we introduce an m.c.n.\ space $H(\Phi,\Psi)$ whose 
matrix norm structure is given by identifying
$M_{n}(H(\Phi,\Psi))$ with $(H\otimes M_{n}, \|\cdot\|_{\Phi,\Psi})$.
We write $H(\Phi)=H(\Phi,\Phi)$ for simplicity.
Before proving that $H(\Phi,\Psi)$ is a homogeneous m.c.n.\ space, we 
prove a useful formula. We denote by $F(K)$ and $U(K)$ 
the subsets of $B(K)$ consisting of all finite-rank operators and all
unitary operators, respectively. If $S$ is a subset of $B(K)$, we denote by $S_+$
the subset of $S$ consisting of positive elements in $B(K)$.

\begin{lem}\label{positive-finite-attainity}
For any operator $T$ we have the equality 
\[ \| T\|^{2}_{\Phi,\Psi}=
\sup\left\{\tr(a\rho_{T} (b))\right\}=\| T^{\ast}\|^{2}_{\Psi^{\ast},\Phi^{\ast}}, \]
where the supremum is taken over all $a,b\in F(K)_+$ 
with $\| a\|_{\Psi^{\ast}}\leq 1$ and $\| b\|_{\Phi}\leq 1$.
\end{lem}

\begin{proof}
Note first that for any $b\in\mathfrak{S}_{\Phi}$ it holds that 
\[\| b\|_{\Phi}=\sup_{\begin{smallmatrix}
    a\in F(K)\\
    \| a\|_{\Phi^{\ast}}\leq 1
    \end{smallmatrix}}|\tr(ab)|, \]
and if $a$ is positive we can choose $b$ to be also 
positive~\cite[proof of Theorem 12.2]{GK}.
The trace duality implies
\[
\|\rho_{T} :\mathfrak{S}_\Phi\rightarrow \mathfrak{S}_\Psi \| 
= \sup_{\| b\|_\Phi\leq 1}\|\rho_{T}(b)\|_{\Psi}
    = \sup_{
        \begin{smallmatrix}
        \| b\|_{\Phi}\leq 1 \\
        \| a\|_{\Psi^{\ast}}\leq 1
        \end{smallmatrix}
    }|\tr(a\rho_{T}(b))|.
\]
If we let $a=u|a|$ and $b=v|b|$ be the polar decompositions 
of $a$ and $b$, respectively, by the Schwarz inequality we have
\begin{eqnarray*}
|\tr(a\rho_{T}(b))| &\leq& 
\tr\left(\sum_{i}\left||a|^{\frac{1}{2}}T_{i}v|b|^{\frac{1}{2}}\right|^{2}\right)^{1/2}
\tr\left(\sum_{i}\left||a|^{\frac{1}{2}}u^{\ast}T_{i}|b|^{\frac{1}{2}}\right|^{2}\right)^{1/2}\\
    &=& \tr(|a|\rho_{T}(v|b|v^{\ast}))^{1/2}
\tr(u|a|u^{\ast}\rho_{T}(|b|))^{1/2}\\
    &\leq& \sup_{\begin{smallmatrix}
        x,y\geq 0 \\
        \| x\|_{\Psi^{\ast}},\| y\|_{\Phi}\leq 1
        \end{smallmatrix}}\tr(x\rho_{T}(y)).
\end{eqnarray*}
Thus
\begin{eqnarray*}
\|\rho_{T} :\mathfrak{S}_\Phi\rightarrow \mathfrak{S}_\Psi \| &=&
\sup_{\begin{smallmatrix}
        x,y\geq 0 \\
        \| x\|_{\Psi^{\ast}},\| y\|_{\Phi}\leq 1
        \end{smallmatrix}}\tr(x\rho_{T}(y))
 = \sup_{\begin{smallmatrix}
        y\geq 0 \\
        \| y\|_{\Phi}\leq 1
        \end{smallmatrix}}\|\rho_{T}(y)\|_{\Psi}\\
 &=& \sup_{\begin{smallmatrix}
        x\in F(K)_+,~y\geq 0 \\
        \| x\|_{\Psi^{\ast}},\| y\|_{\Phi}\leq 1
        \end{smallmatrix}}\tr(x\rho_{T}(y))
 = \sup_{\begin{smallmatrix}
        x\in F(K)_+ \\
        \| x\|_{\Psi^{\ast}}\leq 1
        \end{smallmatrix}}\|\rho_{T^{\ast}}(x)\|_{\Phi^{\ast}}\\
 &=& \sup_{\begin{smallmatrix}
        x,y\in F(K)_+ \\
        \| x\|_{\Psi^{\ast}},\| y\|_{\Phi}\leq 1
        \end{smallmatrix}}\tr(x\rho_{T}(y)).
\end{eqnarray*}
$\Box$
\end{proof}

\begin{prop}\label{M2-satisfied}
The space $H(\Phi,\Psi)$ is an m.c.n.\ space and satisfies 
the Ruan's axiom (M2).
\end{prop}

\begin{proof}
Let $T$ and $S$ be finite sums defined by
\[
T = \sum_{i}\xi_{i}\otimes T_{i},~
S = \sum_{i}\xi_{i}\otimes S_{i},
\]
and let $a,b\in F(K)_{+}$. Then
\begin{eqnarray*}
&& \tr(a\rho_{T+S}(b))\\ 
&=& \sum_{i}\tr\left( a(T_{i}+S_{i})b(T_{i}^{\ast}+S_{i}^{\ast})\right)\\
    &=& \tr (a\rho_{T}(b))+\tr (a\rho_{S}(b))+\sum_{i} (\tr
(a T_{i}b S_{i}^{\ast})+\tr (a S_{i}b T_{i}^{\ast}))\\
    &\leq& \tr (a\rho_{T}(b))+\tr (a\rho_{S}(b))+2\sqrt{\sum_{i}\tr (a T_{i}b T_{i}^{\ast})}
\sqrt{\sum_{i}\tr (a S_{i}b S_{i}^{\ast})}\\
        &=& \tr (a\rho_{T}(b))+\tr (a\rho_{S}(b))+2
\sqrt{\tr(a\rho_{T}(b))\tr(a\rho_{S}(b))}\\
        &=& \left(\tr (a\rho_{T}(b))^{1/2}+\tr (a\rho_{S}(b))
^{1/2}\right)^{2}.
\end{eqnarray*}
Thus $\| T+S\|_{\Phi,\Psi}\leq\| T\|_{\Phi,\Psi}+\| S\|_{\Phi,\Psi}$.
If $T=\xi\otimes A$ is a simple tensor product with $\|\xi\|=1$, then
\[
\|\rho_{T}(x)\|_{\Psi} = \| A x A^{\ast}\|_{\Psi}
        \leq \| A\|\| x\|_{\Psi}\| A\|
        \leq \| A\|\| x\|_{\Phi}\| A\|.
\]
Conversely,
\[ \| T\|_{\Phi,\Psi}^{2} \geq \sup_{p}\| A p A^{\ast}\|_{\Psi}
 = \sup_{p} \| p A^{\ast}A p\|_{\Psi}= \| A\|^{2}, \]
where $p$ runs over all rank one projections.
Thus $\|\xi\otimes A\|_{\Phi,\Psi}=\|\xi\|\|A\|$ 
and hence $H(\Phi,\Psi)$ is an m.c.n.\ space.
Finally, if $X$ and $Y$ are scalar matrices, then
\begin{eqnarray*}
\| X T Y\|_{\Phi,\Psi}^{2} &=& \sup_{a,b}
\frac{\left|\tr\left(\sum_{i}X T_{i}Y a Y^{\ast}T_{i}^{\ast}X^{\ast}b\right)\right|}
{\| a\|_{\Phi}\| b\|_{\Psi^{\ast}}}\\
        &=& \sup_{a,b}\frac{\left|\tr\left(\sum_{i}X T_{i}Y a Y^{\ast}T_{i}^{\ast}X^{\ast}b\right)\right|}
{\| Y a Y^{\ast}\|_{\Phi}\| X^{\ast}b X\|_{\Psi^{\ast}}}\frac{\| Y a Y^{\ast}\|_
{\Phi}\| X^{\ast}b X\|_{\Psi^{\ast}}}{\| a\|_{\Phi}\| b\|_{\Psi^{\ast}}}\\
        &\leq&\| T\|_{\Phi,\Psi}^{2}\| X\|^{2}\|Y\|^{2}.
\end{eqnarray*}
This shows that $H(\Phi,\Psi)$ satisfies Ruan's axiom (M2).
\end{proof}

\begin{lem}\label{homogeneity}
The space $H(\Phi,\Psi)$ is homogeneous.
\end{lem}

\begin{proof}
Let $A\in B(H)$. It suffices to show that for any finite sequence
\[ T=\sum_{i=1}^{m} 
\xi_{i}\otimes T_{i}\in H\otimes M_{n} \]
and $x\in M_{n,+}$, the norm inequality
\[\|\rho_{(A\otimes I)T}(x)\|_{\Psi}\leq\| A\|^{2}\|\rho_{T}(x)\|_{\Psi}.\]
holds.
Let $H_0$ be the finite-dimensional subspace of $H$ spanned by 
$\{A\xi_{i}\}_{i=1}^{m}$ and $\{\eta_{j}\}_{j=1}^{k}$ be an orthonormal basis of $H_0$. 
Then $k\leq m$ and there is an $m\times k$-matrix $B=(b_{i j})$ such that 
$\| B\|\leq\| A\|$ and
$A\xi_{i}=\sum_{j=1}^k b_{i j}\eta_{j}$.
Note that 
\[(A\otimes I_{n})T=\sum_{i}A\xi_{i}\otimes T_{i}=\sum_{j}\eta_{j}\otimes
\left(\sum_{i}b_{i j}T_{i}\right).\]
Thus if we let 
$S_{j}=\sum_{i}b_{i j}T_{i}$ for $1\leq j\leq k$, then
\begin{eqnarray*}
 &&\big\|\rho_{(A\otimes I)T}(x)\big\|_{\Psi}\\ &=& 
    \left\|\sum_{j}S_{j}x S_{j}^{\ast}\right\|_{\Psi}\\
    &=& \left\|\begin{pmatrix}
        S_{1} & \ldots & S_{k} \\
        & & \\
        & \bigcirc &
    \end{pmatrix}
    (I_{k}\otimes x)\begin{pmatrix}
        S_{1}^{\ast} & & \\
        \vdots & & \bigcirc \\
        S_{k}^{\ast} & &
    \end{pmatrix}\right\|_{\Psi}
\end{eqnarray*}
\begin{eqnarray*}
    &=& \left\|(I_{k}\otimes x^{\frac{1}{2}})
    \begin{pmatrix}
        S_{1}^{\ast} & & \\
        \vdots & & \bigcirc \\
        S_{k}^{\ast} & &
    \end{pmatrix}\begin{pmatrix}
        S_{1} & \ldots & S_{k} \\
        & & \\
        & \bigcirc &
    \end{pmatrix}
    (I_{k}\otimes x^{\frac{1}{2}})\right\|_{\Psi} \\
 &=& \left\| (I_{k}\otimes x^{\frac{1}{2}})
        (B^{\ast}\otimes I_{n})\begin{pmatrix}
        T_{1}^{\ast} & & \\
        \vdots & & \bigcirc \\
        T_{m}^{\ast} & &
    \end{pmatrix}\begin{pmatrix}
        T_{1} & \ldots & T_{m} \\
        & & \\
        & \bigcirc &
    \end{pmatrix}(B\otimes I_{n})
    (I_{k}\otimes x^{\frac{1}{2}})\right\|_{\Psi} \\
 & \leq & \| B\|^{2}\left\| (I_{m}\otimes x^{\frac{1}{2}})
    \begin{pmatrix}
        T_{1}^{\ast} & & \\
        \vdots & & \bigcirc\\
        T_{m}^{\ast} & &
    \end{pmatrix}\begin{pmatrix}
        T_{1} & \ldots & T_{m}\\
        & & \\
        & \bigcirc &
    \end{pmatrix}(I_{m}\otimes x^{\frac{1}{2}})
    \right\|_{\Psi}\\
 & \leq & \| A\|^{2}\|\rho_{T}(x)\|_{\Psi}.
\end{eqnarray*}
\end{proof}

Let us see some examples.
Thanks to \cite[Theorem 8.4]{P2}, we have $H(\Phi_{\infty})=R$ 
and $H(\Phi_{1})=C$.

Let $H_{1}$ be a homogeneous Hilbertian m.c.n.\ space and $\Phi$ be an s.n.\ function. 
Mathes and Paulsen \cite[p.1764]{MP} define a new m.c.n.\ space $H_{1,\Phi}$ whose matrix norm is defined by
\[\| T\|_{H_{1},\Phi}=\sup_{x\in \mathfrak{S}_{\Phi},~\|x\|_{\Phi}\leq 1}\| 
(x\otimes I)T\|_{H_{1}},~T\in H\otimes B(K).\]
It is easy to see that $H_{1,\Phi}$ is an m.c.n.\ space.
For example, $H_{\Phi_{\infty}}=H$ and $H_{\Phi_{1}}=H_{\min}$~(see \cite[1.3.\ Proposition]{MP}).
If we are given an s.n.\ function $\Phi$, 
let $\tilde{\Phi}$ be the 2-convexification of $\Phi$ defined by
\[\tilde{\Phi}(a_{1}, \ldots ,a_{n}, \ldots)=\Phi(a_{1}^{2}, \ldots ,a_{n}^{2}, \ldots)^
{1/2},~a\in \hat{k}.\]

\begin{lem}\label{basic-characterization}
For any s.n.\ functions $\Phi$ and $\Psi$ with
$\Phi\geq\Psi$, we have the completely 
isometric identifications
\begin{itemize}
\item $H(\Phi_{1},\Phi)=C_{\widetilde{\Phi^{\ast}}}$,
\item $H(\Phi,\Phi_{\infty})=R_{\tilde{\Phi}}$,
\item $H(\Phi,\Psi)_{\Phi_{2}}=H_{\min}$.
\end{itemize}
In particular, $H(\Phi_{1},\Phi_{\infty})=H_{\min}$.
\end{lem}

\begin{proof}
We first prove the second equation. 
Let $T$ be a finite sum defined by
\[ T=\sum_{i}\xi_{i}\otimes T_{i}\in H\otimes B(K).\]
Then
\[\| T\|_{\Phi,\Phi_{\infty}}^{2}=\sup_{
            \begin{smallmatrix}
             a,b\in F(K)_+\\
            \| a\|_{\Phi},\| b\|_{\Phi_{1}}\leq 1
            \end{smallmatrix}
            }\tr(b\rho_{T}(a)).\]
If we write the spectral decomposition of $b$ by 
$b=\sum_{i}\lambda_{i}p_{i}$ with rank one projections $\{p_{i}\}$, then
\[ \tr (b\rho_{T}(a))=\sum_{i}\displaystyle\lambda_{i}\tr (p_{i}\rho_{T}(a))
 \leq \| b\|_{1}\Max_{i}\{\tr (p_{i}\rho_{T}(a))\}. \]
This shows that $b$ can be replaced by rank one projections. Thus we have
\begin{eqnarray*}
    \| T\|_{\Phi,\Phi_{\infty}}^{2} &=& \sup_a\sup_{
            \begin{smallmatrix}
             p:\mathrm{rank\, one} \\
            \mathrm{projection}
            \end{smallmatrix}
            }\tr(p\rho_{T}(a))\\
        &=& \sup_{p}\|\rho_{T^{\ast}}(p)\|_{\Phi^{\ast}}\\
    &=& \sup_{p}\left\| \begin{pmatrix} T_{1}^{\ast}p & \ldots & T_{n}^{\ast}p \\
        & \bigcirc & \end{pmatrix} 
        \begin{pmatrix} p T_{1} & \\
            \vdots & \bigcirc \\ p T_{n} & 
        \end{pmatrix}\right\|_{\Phi^{\ast}}\\
        &=& \sup_{p}\| (p T_{i}T_{j}^{\ast}p)_{i j}\|_{\Phi^{\ast}}.
\end{eqnarray*}
We write $p$ as $p\zeta = \langle\zeta,\xi\rangle\xi$
with a unit vector 
$\xi\in K$. Then for $\eta=(\eta_{i})_{i=1}^{n}\in K^{n}$ we obtain
\begin{eqnarray*}
(p T_{i}T_{j}^{\ast}p)_{i j}\eta &=& \left(\sum_{j}p T_{i}T_{j}^{\ast}p\eta_{j}\right)_{i}\\
 &=& \left(\sum_{j}\langle\eta_{j},\xi\rangle\langle T_{i}T_{j}^{\ast}\xi,\xi\rangle\xi\right)_{i}\\
 &=& \left(\sum_{j}\langle T_{i}T_{j}^{\ast}\xi,\xi\rangle p\eta_{j}\right)_{i}
 = \big((\langle T_{i}T_{j}^{\ast}\xi,\xi\rangle)_{i j}\otimes p\big)\eta.
\end{eqnarray*}
So it holds that
\[ \| T\|_{\Phi,\Phi_{\infty}}
 = \sup_{\xi}\|(\langle T_{i}T_{j}^{\ast}\xi,\xi\rangle)_{i j}\|_{\Phi^{\ast}} \]
We express any positive operator $a\in\mathfrak{S}_{\Phi}$ with $\| a\|_{\Phi}\leq 1$ in the form 
\[ a=v^{\ast}\diag(a_{1},\ldots,a_{n})v, \]
where $v$ is a unitary matrix and $a_{1}\geq\cdots\geq a_{n}$ are eigenvalues of $a$. 
In the following we denote by $a$ the diagonal matrices $\diag(a_{1},\ldots,a_{n})$.
We write $v=(v(i)_{j})_{i j}$. 
Then $\{v(k)\}_{k=1}^{n}$ is an orthonormal basis of $\mathbb{C}^{n}$.
Thus the above supremum is equal to 
\begin{eqnarray*}
 && \sup_{\xi}\sup_{a\geq 0, \Phi(a)\leq 1}\sup_{v}\left|
\tr(v^{\ast}\mathrm{diag}(a_{1},\ldots ,a_{n})v(\langle T_{i}T_{j}^{\ast}\xi ,\xi\rangle)_{i j})\right| \\
 &=& \sup_{\xi}\sup_{a\geq 0, \Phi(a)\leq 1}\sup_
{\{v(k)\}_{k=1}^{n}}\left|\sum_{k,i,j} a_{k}v(k)_{i}v(k)_{j}
^{\ast}\langle  T_{i}T_{j}^{\ast}\xi,\xi\rangle\right|\\
 &=& \sup_{\xi}\sup_{a\geq 0, \Phi(a)\leq 1}\sup_
{\{v(k)\}_{k=1}^{n}}\left|\left\langle \sum_{k}a_{k}T(v(k))T(v(k))^{\ast}\xi,
\xi\right\rangle\right|,
\end{eqnarray*}
where $T(v(k))$ is defined by
$T(v(k))=\sum_{k=1}^{n} v(k)_{i}T_{i}$.
Hence
\[ \| T\|_{\Phi,\Phi_{\infty}}^{2} = \sup_{a,~\{v(k)\}}\left\|\sum_{k}a_{k}T(v(k))
T(v(k))^{\ast}\right\|=\| T\|_{R_{\tilde{\Phi}}}^{2}. \]

The second equality follows from
\begin{center}
$\| T\|_{R_{\tilde{\Phi}}}=\| T^{\ast}\|_{C_{\widetilde{\Phi^{\ast}}}}$.
\end{center}

The third equality holds since
\begin{eqnarray*}
    \| T\|_{H(\Phi,\Psi)_{\Phi_{2}}} &=& \sup_{
\begin{smallmatrix}
a_{1}\geq\ldots\geq a_{n}\geq 0,~\sum\limits_{i}a_{i}\leq 1 \\
\{v(k)\}_{k=1}^{n} \\
\| x\|_{\Phi}\leq 1,~\| y\|_{\Psi^{\ast}}\leq 1
\end{smallmatrix}
}\left|\tr\left(\sum_{k}a_{k}y T(v(k))x T(v(k))^
{\ast}\right)\right|^{1/2} \\
    &\leq& \sup_{k}\sup_{
\begin{smallmatrix}
\{v(k)\}_{k=1}^{n} \\
\| x\|_{\Phi}\leq 1,~\| y\|_{\Psi^{\ast}}\leq 1
\end{smallmatrix}
}\left|\tr\left( y T(v(k))x T(v(k))^{\ast}\right)\right|^{1/2} \\
    &=& \sup_{v\in \ell_{2}^{n},~\| v\|\leq 1}
\left\|\xi\otimes\left(\sum_{i}v_{i}T_{i}\right)\right\|_{H(\Phi,\Psi)} \\
    &=& \sup_{v\in \ell_{2}^{n},~\| v\|\leq 1}\left\|\sum_{i}v_{i}T_{i}\right\|
=\| T\|_{\min}.
\end{eqnarray*}
Finally, these equalities imply that $H(\Phi_{1},\Phi_{\infty})=C_{\Phi_{2}}=H_{\min}$.
\end{proof}

To check whether $H(\Phi,\Psi)$ is an operator space, it suffices to check whether 
$H(\Phi,\Psi)$ satisfies Ruan's axiom (M1)$^{\prime}$.
The three m.c.n.\ spaces in Lemma \ref{basic-characterization} are clearly operator spaces. 
But not every $H(\Phi ,\Psi)$ is an operator space. 
We give a necessary condition for $H(\Phi,\Psi)$ to be an operator space.

\begin{thm}\label{os-cross}
Let $\Phi$ and $\Psi$ be s.n.\ functions with $\Phi\geq\Psi$.
If the m.c.n.\ space $H(\Phi,\Psi)$ is an operator space, 
then for any $x,y,z\in\mathfrak{S}_{\Phi}$ the following inequality
\[ \frac{\| x\otimes y\|_{\Psi}}{\| x\|_{\Psi}}
\leq\frac{\| z\otimes y\|_{\Phi}}{\| z\|_{\Phi}} \]
holds.
In particular, if $H(\Phi)$ is an operator space, then $\Phi$ is a cross norm.
\end{thm}

\begin{proof}
We may suppose that $x$, $y$, and $z$ are positive diagonal matrices 
in $M_{n}~(n\in\mathbb{N})$ written by
$x=\diag(x_{i})$, $y=\diag(y_{i})$, and $z=\diag(z_{i})$.
For each positive diagonal matrix $w_{i}=\diag(w_{i})\in M_{n}$, 
let $T=\sum_{i, j=1}^{n}\xi_{i}\otimes z_{i}^{1/2}w_{j}^{1/2}e_{i j}$,
Then $\rho_{T}(x)=\sum_{i,j}z_{i}w_{j}x_{j}e_{i i}$ and thus
\[ \| \rho_{T}\|=\sup_{x}\frac{|\tr (x w)|\| z\|_{\Psi}}{\| x\|_{\Phi}}
=\| w\|_{\Phi^{*}}\| z\|_{\Psi}. \]
Let $S$ be the $n$-tuple of $T$. 
Since $\| \rho_{T}\|\geq\| \rho_{S}(x\otimes y)\|_{\Psi}/\|x\otimes y\|_{\Phi}$, 
\[ \| w\|_{\Phi^{*}}\| z\|_{\Psi}\geq
 \frac{|\tr (x w)|\| z\otimes y\|_{\Psi}}{\| x\otimes y\|_{\Phi}}. \]
Taking the supremum over $w$, we obtain the required inequality.
When $\Phi=\Psi$, if we take $x$ or $z$ a rank one projection, 
then we see that $\Phi$ must be a cross norm.
\end{proof}

\begin{ques}
Is the converse of Theorem \ref{os-cross} true?
Namely, if two s.n.\ functions $\Phi$ and $\Psi$ satisfy the conclusion of Theorem \ref{os-cross},
is $H(\Phi,\Psi)$ always an operator space?
\end{ques}

Theorem \ref{os-cross} shows that $H(\Phi)$ is an operator space only if
$\|\cdot\|$ is a cross norm. 
Indeed, we show in Theorem \ref{os-class} that $H(\Phi)$ is an operator space if and only if 
$\Phi$ is the Schatten $p$-norm for some $p\in [1,\infty]$.

\begin{rem}
Let $C_{q}~(1\leq q\leq\infty)$ be the operator space defined by $C_{q}=(C,R)_{1/q}$, and 
we define the operator space 
$S_{p}(C_{q})=(\mathfrak{S}_{1}\hat{\otimes}C_{q},\mathfrak{S}_{\infty}\mt C_{q})_{1/p}$,
where $\hat{\otimes}$ means the operator space projective tensor product (cf. \cite[Section 7]{ER}).
Q. Xu showed in~\cite[Theorem 1]{Xu1} that if we define 
$2\leq p\leq\infty$, $0<\theta <1$, $r,r_{0}(\theta), r_{1}(\theta)$, and $q$ by 
\begin{equation*}
\frac{1}{r}=1-\frac{2}{p},~\frac{1}{r_{0}(\theta)}=\frac{\theta}{2r},~
\frac{1}{r_{1}(\theta)}=\frac{1-\theta}{2r},~\frac{1}{q}=\frac{1-\theta}{p}+\frac{\theta}{p^{\prime}}
\end{equation*}
where $1=1/p+1/p^{\prime}$, then for any $x=(x_{1},x_{2},\ldots ,x_{n})\in \mathfrak{S}_{p}^{n}$,
\begin{equation*}
\| x\|_{S_{p}(C_{q})}=\sup\left\{ \left(\sum_{k}\| a x_{k}b \|_{2}^{2}\right)^{1/2}\right\},
\end{equation*}
where the supremum is taken over all $a\in\mathfrak{S}_{r_{0}(\theta)}$ and $b\in\mathfrak{S}_{r_{1}(\theta)}$ 
with norm one. This is an analogue of $H(\Phi_{p_{1}},\Phi_{q_{1}})$,
where $1/p_{1}=(1-\theta)(1-2/p)$ and $1/q_{1}=1-\theta(1-2/p)$. 
In this case we have $p_{1}\geq q_{1}$.
\end{rem}

\begin{rem}
We can introduce another construction of m.c.n.\ spaces. For any finite sum
$T=\sum_{i}\xi_{i}\otimes T_{i}\in H\otimes M_{n}$ we define
\[\| T\|_{\Phi,\Psi}^{\infty}=\|\rho_{T\otimes I_{\infty}}\colon\mathfrak{S}_{\Phi}
\rightarrow\mathfrak{S}_{\Psi}\|^{1/2},\]
where $T\otimes I_{\infty}$ acts on $B(K\otimes \ell_{2})$ and $K\otimes \ell_{2}$ is 
identified with a separable infinite-dimensional Hilbert space. 
Then we denote by $H^{\sharp}(\Phi,\Psi)$ the m.c.n.\ space
whose matrix norm structure is given by the family
$(H\otimes M_{n},\|\cdot\|_{\Phi,\Psi}^{\infty})$.
There is a case where $H^{\sharp}(\Phi,\Psi)$ is an operator space though 
$H(\Phi,\Psi)$ is not an operator space. 
Let $\Phi$ be the KyFan 2-norm, 
that is, $\Phi(a)=a_{1}^{\ast}+a_{2}^{\ast}$. 
Then $H(\Phi)$ is not an operator space. 
Indeed, for $x=\diag(1,1)\in M_{2}$ it holds that $\| x\otimes x\|_{\Phi}=2$,
but $\| x\|_{\Phi}^{2}=4$.
To determine $H^{\sharp}(\Phi)$, if we are given Hilbertian operator 
spaces $H_{1}$ and $H_{2}$ with the common base space $H$,
we define the matricially normed
space $H_{1}\bigvee H_{2}$ with the base space $H$ by
\[\| x\|_{M_{n}(H_{1}\bigvee H_{2})}=\Max\{\| x\|_{M_{n}(H_{1})},\| x\|_{M_{n}(H_{2})}\}.\]
It is easy to see that $H_{1}\bigvee H_{2}$ is an operator space.
\end{rem}

\begin{prop}\label{ex01}
Let $\Phi$ be an s.n.\ function defined by 
\[\Phi(a)=a_{1}^{\ast}+\theta a_{2}^{\ast}~(0<\theta\leq 1).\]
Then $H^{\sharp}(\Phi)$ is an operator space equal to 
$H^{\sharp}(\Phi)=H(\Phi_{\infty})\bigvee H(\Phi_{1},\Phi)$.
\end{prop}

\begin{proof}
Let $T$ be a finite sum defined by $T=\sum_{i}\xi_{i}\otimes T_{i}$.
For any $x\in F(K)_+$ we write its spectral decomposition as
$x=\sum_{j=1}^{m}s_{j}(x)p_{j}$.
Then if we let 
\[y=s_{1}(x)p_{1}+s_{2}(x)\sum_{j=2}^{m}p_{j},\]
then $y$ satisfies $\| y\|_{\Phi}=\| x\|_{\Phi}$ and $x\leq y$. Thus we have
\begin{eqnarray*}
\|\rho_{T\otimes I_{\infty}}\|_{\Phi} &=& \sup_{\frac{1}{1+\theta}\leq\alpha\leq 
1}\sup_{p,q}\left\|\rho_{T\otimes I_{\infty}}(\alpha
p+\frac{1-\alpha}{\theta}q)\right\|_{\Phi}\\
    &=& \sup_{p,q}\Max\left\{\|\rho_{T\otimes I_{\infty}}(p)\|_{\Phi},\frac
{\|\rho_{T\otimes I_{\infty}}(p+q)\|_{\Phi}}{1+\theta}\right\},
\end{eqnarray*}
where $p$ runs over all rank one projections and $q$ runs over all finite rank projections 
orthogonal to $p$.
Now for fixed $p$, it is clear that
\[\|\rho_{T\otimes I_{\infty}}(p+q)\|_{\Phi}\leq
    (1+\theta)\left\|\sum_{i}T_{i}T_{i}^{\ast}\right\|\]
for any projection $q$ orthogonal to $p$.
To show the converse, represent $p$ as $p\eta = \langle\eta,\xi\rangle\xi$ 
with a unit vector $\xi$ and write 
\[\xi=\sum_{i=1}^{n}\phi_{i}\otimes\psi_{i},~\phi_{i}\in \ell_{2}^{n},~\psi_{i}\in \ell_{2}.\]
If we take a projection $r\in B(\ell_{2})$ such that the rank of $r$ is not less 
than 2 and orthogonal to the vectors $\{\psi_{i}\}$ and let $q=I_{n}\otimes r$, then 
we have
\[ \|\rho_{T\otimes I_{\infty}}(p+q)\|_{\Phi}
\geq \left\|\sum_{i}T_{i}T_{i}^{\ast}\otimes r\right\|_{\Phi}
= (1+\theta)\left\|\sum_{i}T_{i}T_{i}^{\ast}\right\|. \]
Thus
\[\sup_{p,q}\Max\left\{\|\rho_{T\otimes I_{\infty}}(p)\|_{\Phi},\frac
{\|\rho_{T\otimes I_{\infty}}(p+q)\|_{\Phi}}{1+\theta}\right\}
 =\Max\left\{\| T\|_{\Phi_{1},\Phi}^{2},\left\|\sum_{i}T_{i}T_{i}^{\ast}\right\|\right\}.\]
\end{proof}

\begin{ques}
Is $H^{\sharp}(\Phi,\Psi)$ always an operator space?
\end{ques}

As we see below, 
for many two distinct s.n.\ functions $\Phi\neq\Psi$, 
the m.c.n.\ space $H(\Phi,\Psi)$ is an operator space.
Pisier~\cite[Theorem 8.4]{P2} showed the completely isometrically isomorphism
$H(\Phi_{p},\Phi_{p})=R(\theta)$, where $1\leq p\leq \infty$ and $\theta=p^{-1}$.
We consider whether $H(\Phi_{p},\Phi_{q})$ is an operator space 
for general $p$ and $q$ with $1\leq p\leq q\leq \infty$. In the case of $p=1$ 
or $q=\infty$, $H(\Phi_{p},\Phi_{q})$ is an operator space from Lemma \ref{basic-characterization}.
To deal with the case $1\leq p\leq 2\leq q\leq\infty$ we need the following notion.

\begin{defn}
Let $\Phi$ be an s.n.\ function. We call $\Phi$ a Q-norm if there 
is an s.n.\ function $\Upsilon$ such that $\tilde{\Upsilon}=\Phi$, 
and $\Phi$ is a Q$^{\ast}$-norm if $\Phi$ is an adjoint of some Q-norm.
In other words, an s.n.\ function $\Phi$ is a Q-norm if 
there is an s.n.\ function $\Upsilon$ such that for any $A\in\mathfrak{S}_{\Phi}$, 
the norm equality
\[ \| A\|_{\Phi}^{2}=\| A^{\ast}A\|_{\Upsilon}\]
is satisfied. 
Note that a Q-norm is smaller than or equal to the Schatten 2-norm and 
a Q$^{\ast}$-norm is greater than or equal to the Schatten 2-norm.
For example, the Schatten $p$-norm $\Phi_{p}$ is a Q-norm when 
$2\leq p\leq\infty$ and is a Q$^{\ast}$-norm when $1\leq p\leq 2$.
The Lorentz ideal $\Phi_{p,q}$ is a Q-norm if $2\leq q$.
We use the following lemma.
\end{defn}

\begin{lem}\label{matrix-partitioned-norm}\cite[Proposition 3]{BK}
Let $\Phi$ be a Q$^{\ast}$-norm and $y=\bigl( \begin{smallmatrix}
    y_{1} & y_{2}\\
    y_3 & y_4
    \end{smallmatrix}\bigr)$ with $y_{i}\in M_{n}~(i=1,2,3,4$ and $n\in \mathbb{N})$. 
Then we have the inequality
\[\sum_{i=1}^4\| y_{i}\|_{\Phi}^{2}\leq \| y\|_{\Phi}^{2}.\]
\end{lem}

\begin{thm}
Let $\Phi$ be a Q$^{\ast}$-norm and $\Psi$ be a Q-norm. 
Then $H(\Phi,\Psi)$ is an 
operator space.
\end{thm}

\begin{proof}
It suffices to check the Ruan's axiom (M1)$^{\prime}$. Let
$T$ and $S$ be finite sums given by
\[ T=\sum_{i=1}^{k} \xi_{i}\otimes T_{i}\in M_{m}(H(\Phi,\Psi))~\text{and}~
S=\sum_{i=1}^{l} \xi_{i}\otimes S_{i}\in M_{n}(H(\Phi,\Psi)). \]
Since for any $t\in\mathbb{N}$ it follows that 
$\| T\oplus 0_{t}\|_{H(\Phi,\Psi)}=\| T\|_{H(\Phi,\Psi)}$,
we may assume that $m=n$ and clearly that $k=l$. 
Take matrices $y$ and $z$ given by
\[ y=
\left( \begin{matrix}
    y_{1} & y_{2}\\
    y_{3} & y_{4}
    \end{matrix}\right),~z=\left(\begin{matrix}
    z_{1} & z_{2}\\
    z_{3} & z_{4}
    \end{matrix}\right)\in M_{2n,+} \]
with $y_{j},z_{j}\in M_{n}(i=1,2,3,4)$. Then we have
\begin{eqnarray*}
    && \left|\tr\left(\sum_{i} \left(\begin{matrix}
    T_{i} & 0\\
    0 & S_{i}
    \end{matrix}\right) y\left(\begin{matrix}
    T_{i}^{\ast} & 0\\
    0 & S_{i}^{\ast}
    \end{matrix}\right) z\right)\right| \\
    &=& \left|\sum_{i}\tr\left(T_{i}y_{1}T_{i}^{\ast}z_{1}+T_{i}y_{2}S_{i}^{\ast}
z_{3}+S_{i}y_{3}T_{i}^{\ast}z_{2}+S_{i}y_{4}S_{i}^{\ast}z_{4}\right)\right|\\
    &\leq& \mathrm{Max}\left\{\| T\|_{H(\Phi,\Psi)}^{2}, \| S\|_{H(\Phi,\Psi)}^{2}\right\}
\sum_{j=1}^{4}\| y_{j}\|_{\Phi}\| z_{j}\|_{\Psi^{\ast}}\\
    &\leq& \mathrm{Max}\left\{\| T\|_{H(\Phi,\Psi)}^{2}, \| S\|_{H(\Phi,\Psi)}^{2}\right\}
\left\{\sum_{j=1}^{4}\| y_{j}\|_{\Phi}^{2}\right\}^{1/2}\left\{\sum_{j=1}^{4}\| z_{j}\|_{\Psi^{\ast}}^{2}\right\}^{1/2}\\
    &\leq& \mathrm{Max}\left\{\| T\|_{H(\Phi,\Psi)}^{2}, \| S\|_{H(\Phi,\Psi)}^{2}\right\}
\| y\|_{\Phi}\| z\|_{\Psi^{\ast}}.
\end{eqnarray*}
In the third line we use the Schwarz inequality \cite[Theorem IX.5.11]{Ba} and
in the last line we do the preceding lemma. This shows that (M1)$^{\prime}$ holds.
\end{proof}


\section{Completely bounded mappings between $H(\Phi, \Psi)$s.}\label{cbm}

We consider the relationship between the m.c.n.\ spaces $H(\Phi,\Psi)$ and the 
space of completely bounded mappings between them.
It is possible to describe the space $CB(H(\Phi_{\infty}),H(\Phi,\Psi))$ in 
terms of the multiplicator norm, which was discussed by \cite{AMS} in the 
case of rearrangement invariant spaces on the interval $[0,1]$.

\begin{thm}\label{cb-from-infty}
Let $\Phi ,\Psi$ be s.n.\ functions with $\Phi\geq\Psi$
and $x\in B(H)$. Then
\[\| x\|_{CB(R,H(\Phi ,\Psi))}=
\left(\sup_{a\in\mathfrak{S}_{\Phi}}\frac{\big\| |x|^{2}\otimes 
a\big\|_{\Psi}}{\| a\|_{\Phi}}\right)^{1/2}.\]
In particular, if $\Phi$ and $\Psi$ satisfy the condition of Theorem \ref{os-cross},
then we have the isometric isomorphisms
$CB(R,H(\Phi ,\Psi))=\mathfrak{S}_{\tilde{\Psi}}$ and
$CB(C,H(\Phi ,\Psi))=\mathfrak{S}_{\tilde{\Phi^{\ast}}}$.
\end{thm}

\begin{proof}
Let $x=\diag(\lambda_{1},\ldots ,\lambda_{n}),~\lambda_{1}\geq\ldots 
\geq\lambda_{n}\geq 0$ be a positive diagonal matrix. Then from the definition
\[ \| x\|_{CB(R,H(\Phi ,\Psi))}=
\sup_{T\in R,~a\in\mathfrak{S}_{\Phi,+},~\| a\|_{\Phi}\leq 1}
\left\{\left\|\sum_{i}\lambda_{i}^{2}T_{i}a T_{i}^{\ast}\right\|_{\Psi}
^{1/2}\right\}. \]
If $\| T\|_{R}\leq 1$, then
$\left\|\sum_{i}T_{i}T_{i}^{\ast}\right\|\leq 1$
and thus it follows that $(T_{i}^{\ast}T_{j})_{i j}\leq I$. Hence we have
\begin{eqnarray*}
\left\|\sum_{i=1}^{n}\lambda_{i}^{2}T_{i}a T_{i}^{\ast}\right\|_{\Psi}
 &=& \left\|
    \begin{pmatrix}
     T_{1} & \ldots & T_{n} \\
     & & \\
     & \bigcirc & 
    \end{pmatrix}
    \diag(\lambda_{1}^{2}a,\ldots,\lambda_{n}^{2}a)
    \begin{pmatrix}
        T_{1}^{\ast} & &\\
    \vdots & \bigcirc &\\
    T_{1}^{\ast} & &
    \end{pmatrix}
    \right\|_{\Psi}\\
 &=& \left\|
    \diag(\lambda_{1}a^{\frac{1}{2}},\ldots,\lambda_{n}a^{\frac{1}{2}})
    (T_{i}^{\ast}T_{j})
    \diag(\lambda_{1}a^{\frac{1}{2}},\ldots,\lambda_{n}a^{\frac{1}{2}})
    \right\|_{\Psi}\\
 &\leq& \| |x|^{2}\otimes a\|_{\Psi}.
\end{eqnarray*}
To show the converse, take a family $\{ T_{i}\}_{i=1}^{n}$ such that $T_{i}^{\ast}T_{j}
=\delta_{i j}I$, where $\delta_{i j}$ is the Kronecker delta.

When $\Phi$ and $\Psi$ satisfy the condition of Theorem \ref{os-cross}, we have
\[ \||x|^{2}\otimes a\|_{\Psi}\leq \||x|^{2}\|_{\Psi}\| a\|_{\Phi}=\| x\|_{\tilde{\Psi}}^{2}\| a\|_{\Phi} \]
and thus $\| x\|_{CB(H(\Phi_{\infty}),H(\Phi ,\Psi))}\leq \| x\|_{\tilde{\Psi}}$.
The converse is verified by putting $a$ to be any rank one projection.
The last assertion is obtained from Lemma \ref{positive-finite-attainity}.
\end{proof}

Other important Hilbertian operator spaces are $H_{\min}$ and $OH$.
Let us see the space $CB(H_{\min}, H(\Phi_{p},\Phi_{q}))$ next.
When $p=q$, this space can be identified with $\mathfrak{S}_{2}$.

\begin{thm}\label{cb-norm-min-R}
For each $\theta\in [0,1]$, the space $CB(H_{\min},R(\theta))$ 
coincides with $\mathfrak{S}_{2}$ up to equivalence of norm.
\end{thm}
\begin{proof}
Mathes proved this theorem when $\theta=0$ or $1$ (see \cite[Proposition 6]{Ma}).
We use this result and the complex interpolation theory.
Since the space of completely bounded mappings between homogeneous 
m.c.n.\ spaces is an operator ideal, it suffices to check the cb-norm of the 
matrices of the diagonal form $A=\diag(\lambda_{1},\ldots ,\lambda_{n}),~
\lambda_{1}\geq\ldots \geq\lambda_{n}\geq 0$. We denote by $\| A\|_{cb}$ the c.b.\ norm 
of $A\colon H_{\min}\rightarrow R(\theta)$.
First we note that
\[ \| A\|_{cb}=\sup_{T}\frac{\left\|\sum_{i}\xi_{i}\otimes\lambda_{i}T_{i}\right\|_{R(\theta)}}
{\| T\|_{\min}}. \]
Thus by the complex interpolation property it follows that
\[ \| A\|_{cb}\leq\sup_{T}\left\{\frac{\left\|\sum_{i}\xi_{i}\otimes\lambda_{i}T_{i}\right\|_{R}}
{\| T\|_{\min}}\right\}^{1-\theta}\left\{\frac{\left\|\sum_{i}\xi_{i}\otimes\lambda_{i}T_{i}\right\|_{C}}
{\| T\|_{\min}}\right\}^{\theta}\leq\left(\sum_{i}\lambda_{i}^{2}\right)^{1/2}, \]
where we use the case of $\theta=0,1$.
To show the converse inequality, we use the spin system $\{ U_{i}\}$. This system 
is an $n$-tuple of unitary self-adjoint operators such that
\[\forall i\neq j,~U_{i}U_{j}+U_{j}U_{i}=0\]
(cf. \cite[p.76]{P2}). The spin system satisfies
\[\left\|\sum_{i}\eta_{i}U_{i}\right\|\leq\sqrt{2}\left(\sum_{i}|\eta_{i}|^{2}\right)
^{1/2},~\forall (\eta_{i})\in\mathbb{C}^{n}\]
and
\[\left\|\sum_{i}\lambda_{i}^{2}U_{i}\otimes U_{i}\right\|=\left(\sum_{i}\lambda_{i}^{2}\right)
^{1/2}.\]
The first property implies that 
\[\left\|\sum_{i}\xi_{i}\otimes U_{i}\right\|_{\min}\leq\sqrt{2}.\]
The complex interpolation duality leads the isomorphism $R(\theta)^{\ast}=R(1-\theta)$. 
Using this we obtain
\begin{eqnarray*}
\sum_{i}\lambda_{i}^{2} = \left\|\sum_{i}\lambda_{i}^{2}U_{i}\otimes U_{i}\right\|
    &\leq& 2\frac{\left\|\sum_{i}\xi_{i}\otimes\lambda_{i}U_{i}\right\|_{R(\theta)}}{\big\| U\big\|_{\min}}
\frac{\left\|\sum_{i}\xi_{i}\otimes\lambda_{i}U_{i}\right\|_{R(1-\theta)}}{\big\| U\big\|_{\min}}\\
    &\leq& 2\| A\|_{CB(H_{\min},R(\theta))}\| A\|_{CB(H_{\min},R(1-\theta))}\\
    &\leq& 2\| A\|_{CB(H_{\min},R(\theta))}\left(\sum_{i}\lambda_{i}^{2}\right)^{1/2}.
\end{eqnarray*}
Thus $\| A\|_{2}\leq 2\| A\|_{CB(H_{\min},R(\theta))}$.
\end{proof}

To deal with the case $p\neq q$, we need the following lemma.

\begin{lem}\label{interpolation01}
Let $1\leq p\leq q\leq\infty$ and take $\theta ,\psi\in [0,1]$ such that
\[\begin{cases}
    1/p=1-\psi+\theta\psi\\
    1/q=\theta\psi.
    \end{cases}\]
Then for every $T\in H(\Phi_{p},\Phi_{q})$,
\[\| T\|_{\Phi_{p},\Phi_{q}}\leq\| T\|_{(H_{\min}, R(\theta))_{\psi}}.\]
\end{lem}
\begin{proof}
For each $t\in [0,1]$, take positive numbers $p_{t}$ and $q_{t}$ such that
\[ \frac{1}{p_{t}}=1-t+\theta t,~\frac{1}{q_{t}}=\theta t\]
and let $q_{t}'=\left(1-1/q_{t}\right)^{-1}$. 
We define a family of bilinear mappings 
$f_{t}\colon\mathfrak{S}_{2q_{t}'}
\times\mathfrak{S}_{2p_{t}}\rightarrow \ell_{2}(\mathfrak{S}_{2})$ by 
$f_{t}(a,b)=(a T_{i}b)_{i}$ for $0\leq t\leq 1$.
Then Lemma \ref{basic-characterization} shows that $\| f_0\|=\| T\|_{\min}$ and Pisier~\cite[Theorem 8.4]{P2} shows
$\| f_1\|=\| T\|_{R(\theta)}$. Thus the multilinear interpolation 
(see \cite[10.2]{C}) implies that $\| T\|_{\Phi_p,\Phi_q}=\| 
f_{\psi}\|\leq\| T\|_{(H_{\min},R(\theta))_{\psi}}$.
\end{proof}

\begin{thm}
Let $1\leq p\leq q\leq\infty$. We have a contractive embedding of 
$\mathfrak{S}_{r}$ into $CB(H_{\min}, 
H(\Phi_{p},\Phi_{q}))$, where $r=2/(1/q-1/p+1)$.
\end{thm}

\begin{proof}
Let $A=\diag(\lambda_{1},\ldots ,\lambda_{n}),~\lambda_{1}\geq\ldots 
\geq\lambda_{n}\geq 0$. Then,
\begin{eqnarray*}
\| A\|_{CB(H_{\min}, H(\Phi_{p},\Phi_{q}))} &\leq& \| A\|_{CB(H_{\min},(H_{\min}, 
R(\theta))_{\psi})}\\
    &\leq& \| A\|_{(CB(H_{\min},H_{\min}),CB(H_{\min},R(\theta)))_{\psi}}\\
    &\leq& \| A\|_{(\mathfrak{S}_{\infty},\mathfrak{S}_{2})_{\psi}}
    = \| A\|_{r}.
\end{eqnarray*}
In the first step we use Lemma \ref{interpolation01} and in the third we use Theorem \ref{cb-norm-min-R}.
\end{proof}

We observe the c.b.\ norm of the mappings from $OH$ to $H(\Phi_{p},\Phi_{q})$.
\begin{thm}
Let $1\leq p\leq q\leq \infty$. Then
\begin{equation*}
CB(OH,H(\Phi_{p},\Phi_{q}))=\begin{cases}
    \mathfrak{S}_{4(1-2/p)^{-1}} & (p\geq 2) \\
    B(H) & (p\leq 2\leq q)\\
    \mathfrak{S}_{4(2/q-1)^{-1}} & (q\leq 2)
    \end{cases}
\end{equation*}
with equal norms.
\end{thm}

\begin{proof}
The second case is obvious and the third one follows from 
Lemma \ref{positive-finite-attainity} and the first one.
We show the first case. Let $A=\diag(\lambda_{1},\ldots ,\lambda_{n})$ be a 
diagonal operator with $\lambda_{1}\geq\cdots\geq\lambda_{n}\geq 0$.
Xu showed in \cite[Lemma 5.9]{Xu2} that if $1\leq p\neq q \leq \infty$, then 
$CB(H(\Phi_{p},\Phi_{p}),H(\Phi_{q},\Phi_{q}))=\mathfrak{S}_{2 p q/|p-q|}$.
From this result it clearly follows that for any operator $A$,
\begin{equation*}
\| A\|_{CB(OH,H(\Phi_{p},\Phi_{q}))}\leq \| A\|_{CB(OH,H(\Phi_{p},\Phi_{p}))}
=\| A\|_{4(1-2/p)^{-1}}.
\end{equation*}
To show the converse, for a positive diagonal matrix $B=\diag(b_{1},\ldots ,b_{n})$,
let $T_{B,i}=b_{i} e_{1i}\in M_{n}~(i=1,\ldots ,n)$. Then
\begin{align*}
\left\|\sum_{i=1}^{n}\xi_{i}\otimes T_{B,i}\right\|_{OH}^4 
    &= \left\|\sum_{i=1}^{n}T_{B,i}\otimes \bar{T}_{B,i}\right\|_{\min}^{2} \\
    &= \left\|\sum_{i,j=1}^{n}b_{i}^{2}b_{j}^{2}(e_{1i}\otimes e_{1i})
        (e_{j1}\otimes e_{j1})\right\|_{M_{n}\otimes M_{n}} \\
    &= \left\|\sum_{i=1}^{n}b_{i}^{4}e_{11}\otimes e_{11}\right\|_{M_{n}\otimes M_{n}}=\sum_{i=1}^{n}b_{i}^{4}.
\end{align*}
However, if we let $C$ be a positive diagonal matrix 
$\diag(c_{1},\ldots,c_{n})$, then we have
\[ \left\| A\left(\sum_{i=1}^{n}\xi_{i}\otimes T_{B,i}\right)\right\|_{H(\Phi_{p},\Phi_{q})} 
    \geq \sup_{C}\frac{\left|\sum_{i=1}^{n}\lambda_{i}^{2}b_{i}^{2}c_{i}\right|^{1/2}}
            {\left(\sum_{i=1}^{n}c_{i}^{p}\right)^{1/p}}.\]
Taking the supremum for $B$ in the unit ball of $OH$, we obtain
\[ \| A\|_{CB(OH,H(\Phi_{p},\Phi_{q}))}
\geq \sup_{C}\frac{\left|\sum_{i=1}^{n}\lambda_{i}^{4}c_{i}^{2}\right|^{1/4}}
            {\left(\sum_{i=1}^{n}c_{i}^{p}\right)^{1/p}}
=\| A\|_{4(1-2/p)^{-1}}. \]
\end{proof}


\section{Multiplicator in operator ideals}\label{multi}

In this section we show that the m.c.n.\ space $H(\Phi)$ is
an operator space if and only if $\Phi$ is the Schatten norm.

In view of the result of Theorem \ref{cb-from-infty}, for an s.n.\ function $\Phi$ we 
consider the following two conditions
\begin{enumerate}
\item[$(\ast)$]$\exists c_{1}\geq 
0,~\|x\otimes y\|_{\Phi}\leq c_{1}\| x\|_{\Phi}\| y\|_{\Phi}~
\text{for any}~x~\text{and}~y$;
\item[$(\ast\ast)$]$\exists c_{2}\geq 
0,~\|x\otimes y\|_{\Phi}\geq c_{2}\| x\|_{\Phi}\| y\|_{\Phi}~
\text{for any}~x~\text{and}~y$.
\end{enumerate}
Note that if an s.n.\ function $\Phi$ satisfies $(\ast)$, 
its adjoint $\Phi^{\ast}$ satisfies $(\ast\ast)$ for $c_{2}$ with $c_{1}c_{2}=1$.
The Schatten $p$-norm is a cross norm and satisfies both $(\ast)$ 
and $(\ast\ast)$ with $c_{1}=c_{2}=1$.

Let $\Phi$ and $\Psi$ be s.n.\ functions with $\Phi\geq\Psi$ and 
$x\in B(\ell_{2})$ such that 
\[\sup_a\frac{\| x\otimes a\|_{\Psi}}{\| a\|_{\Phi}}<\infty.\] 
We denote by 
$M_{\Phi,\Psi}(x)$ the multiplicator from $\mathfrak{S}_{\Phi}$ to 
$\mathfrak{S}_{\Psi}$ defined by 
\[ M_{\Phi,\Psi}(x)(a)=x\otimes a. \] 
For an s.n.\ function $\Phi$, we denote by $\mathcal{M}(\mathfrak{S}_{\Phi})$
the space consisting of $x\in B(\ell_{2})$ with $M_{\Phi,\Phi}(x)$ is bounded.
We equip $\mathcal{M}(\mathfrak{S}_{\Phi})$ with the norm 
$\| M_{\Phi,\Phi}(x)\|$.
It holds that
\[ \| x\|_{\Psi}=\frac{\| x\otimes e_{11}\|_{\Psi}}{\| e_{11}\|_{\Phi}}
\leq\| M_{\Phi,\Psi}(x)\|. \]
In case of the Schatten norm $(1\leq p\leq q\leq\infty)$, we have
\[ \| M_{\Phi_{p},\Phi_{q}}(x)\|=\| x\|_{q}. \]
If an s.n.\ function $\Phi$ satisfies $(\ast)$, then
\[ \| M_{\Phi,\Phi}(x)\|\leq c_{1}\| x\|_{\Phi}, \]
and thus $\Phi$ satisfies $(\ast)$ if and only if $\| x\|_{\Phi}$ is 
equivalent to $\| M_{\Phi,\Phi}(x)\|$.
Since $M_{\Phi,\Phi}(x)M_{\Phi,\Phi}(y)=M_{\Phi,\Phi}(x\otimes y)$, we have
\[\| M_{\Phi,\Phi}(x\otimes y)\|\leq\| M_{\Phi,\Phi}(x)\|\| M_{\Phi,\Phi}(y)\|.\]
The multiplicator is discussed in \cite{AMS} for the rearrangement invariant 
space on [0,1].

The conditions $(\ast)$ and $(\ast\ast)$ are closely related to the Schatten norm.

\begin{lem}\label{logarithemic-order}
If an s.n.\ ideal $\mathfrak{S}_{\Phi}$ satisfies $(\ast)$ or 
$(\ast\ast)$, then the limit 
\[p=\lim_{n\rightarrow\infty}\frac{\log n}{\log\| P_{n}\|_{\Phi}}\in [1,\infty] \]
exists, where $P_{n}$ stands for any rank $n$ projection.
\end{lem}

\begin{proof}
We prove the statement in the case that $(\ast)$ holds. 
In the case of $(\ast\ast)$ the proof is similar. 
By the hypothesis, for fixed $m\in\mathbb{N}$,
\[\| P_{m^k}\|_{\Phi}\leq c_{1}^{k-1}\|P_{m}\|_{\Phi}^k,\quad\forall k\in\mathbb{N}.\]
If $\{ t_{i}\}_{i=1}^{\infty}$ is a subsequence of $\mathbb{N}$, we can take a 
non-decreasing sequence $\{k_{i}\}_{i=1}^{\infty}$ in $\mathbb{N}$ which tends 
to infinity such that
$m^{k_{i}}\leq t_{i}< m^{k_{i}+1}$.
Thus we have
\[\frac{\log t_{i}}{\log\| P_{t_{i}}\|_{\Phi}}\geq\frac{\log m^{k_{i}}}
{\log\| P_{m^{k_{i}+1}}\|_{\Phi}}\geq
\frac{k_{i}\log m}{k_{i}\log c_{1}+(k_{i}+1)\log\| P_{m}\|_{\Phi}}.\]
Since $\{ t_{i}\}_{i=1}^{\infty}$ is arbitrary, it follows that
\[\liminf_{n\rightarrow\infty}\frac{\log n}{\log\| P_{n}\|_{\Phi}}\geq\frac{\log m}
{c_{1}+\log\| P_{m}\|_{\Phi}}.\]
This implies
\[\liminf_{n\rightarrow\infty}\frac{\log n}{\log\| P_{n}\|_{\Phi}}\geq\limsup_
{m\rightarrow\infty}\frac{\log m}{\log\| P_{m}\|_{\Phi}}\]
and the limit exists.
\end{proof}

\begin{thm}\label{p-Phi-inequality}
Suppose that an s.n.\ ideal $\mathfrak{S}_{\Phi}$ satisfies $(\ast)$ or 
$(\ast\ast)$ and let $p$ be as in the preceding lemma. 
Then the following statements hold.

(1) if $\mathfrak{S}_{\Phi}$ satisfies $(\ast)$, then 
\[\| x\|_{p}\leq c_{1}\| x\|_{\Phi},~\forall x\in\mathfrak{S}_{\Phi}.\]

(2) if $\mathfrak{S}_{\Phi}$ satisfies $(\ast\ast)$, then 
\[c_{2}\| x\|_{\Phi}\leq \| x\|_{p},~\forall x\in\mathfrak{S}_{\Phi}.\]
In particular, if $\Phi$ is a cross norm, then $\Phi=\Phi_{p}$.
\end{thm}

\begin{proof}
Let $x=\diag(\lambda_{1},\ldots,\lambda_{m}),~\lambda_{1}\geq\ldots 
\geq\lambda_{m}\geq 0$ be a diagonal matrix and let 
\[ x^{\otimes n}=\sum_{i=1}^{N}t_{i}e_{i}\]
be the spectral decomposition of the $n$-fold tensor product of $x$.
In the above inequality, $N$ is dominated by $\binom{m+n-1}{m-1}$. 
If we let $p_{j}$ be the $j$-th sum of the $e_{i}$'s given by
$p_{j}=\sum_{i=1}^{j}e_{i}$,
then for all $j$ we have
\[ \sum_{i=1}^{n}t_{i}e_{i}=\sum_{i=1}^{n}(t_{j}-t_{j-1})p_{j}\geq t_{j}p_{j}. \]
Thus it holds that
\[ \Max_{j}\big\{t_{j}\| p_{j}\|_{\Phi}\big\}\leq
\| x^{\otimes n}\|_{\Phi}\leq N\Max_{j}\big\{t_{j}\| p_{j}\|_{\Phi}\big\} \]
and hence
\[ \Max_{j}\{(t_{j}\| p_{j}\|_{\Phi})^{1/n}\}\leq
\| x^{\otimes n}\|_{\Phi}^{1/n}\leq 
N^{1/n}\Max_{j}\{(t_{j}\| p_{j}\|_{\Phi})^{1/n}\}. \]
Note that from the above inequality, if $\Phi=\Phi_{p}$, then
\[\| x\|_{p}=\lim_{n\rightarrow\infty}\Max_{j}
\big\{t_{j}^{1/n}(\rank p_{j})^{1/(p n)}\big\},\]
which proves (1). The proof of (2) is similar. By the preceding lemma, for any 
$\varepsilon\geq 0$, there exists a $D\geq 0$ such that
\[\| p_{j}\|_{\Phi}\geq 
D(\rank p_{j})^{1/(p+\varepsilon)},~\text{for all}~j\in\mathbb{N}.\]
$(\ast)$ implies that $\| x^{\otimes n}\|_{\Phi}\leq c_{1}^{n-1}\| x\|_{\Phi}^{n}$, so 
that
\begin{eqnarray*}
c_{1}\| x\|_{\Phi} &\geq& \| x^{\otimes n}\|_{\Phi}^{1/n}\\
    &\geq& \Max_{j}\{(t_{j}\| p_{j}\|_{\Phi})^{1/n}\}\\
    &\geq& \Max_{j}\big\{(D t_{j})^{1/n}
(\rank p_{j})^{1/\{ (p+\varepsilon)n\} }\big\}.
\end{eqnarray*}
The last term converges to $\| x\|_{p+\varepsilon}$ as $n\rightarrow\infty$.
\end{proof}

From Theorem \ref{p-Phi-inequality} and Theorem \ref{os-cross}, we obtain the following corollary.

\begin{cor}\label{os-class}
Let $\Phi$ be an s.n.\ function. 
The m.c.n.\ space $H(\Phi)$ is an operator space if and only if 
$\Phi$ is some Schatten $p$-norm~$(1\leq p\leq\infty)$.
\end{cor}

\begin{rem}
Let $X$ be a rearrangement invariant function space $X$ on the interval $[0,1]$
(cf. \cite[Section 2]{LT}).
For $s>0$, let $\sigma_s$ be the dilation operator given by
\[\sigma_{s} x(t)=x(t/s)1_{[0,\max\{1,s\}]}~(t\in[0,1],~x\in X).\]
This operator is well defined on $X$ and $\|\sigma_{s}\|\leq\max\{1,s\}$.
The Boyd indices $\alpha_X$ and $\beta_X$ of $X$ are defined by
\[\alpha_X=\lim_{s\rightarrow 0}\frac{\log\|\sigma_s\|_{X\rightarrow X}}{\log s},
\quad \beta_X=\lim_{s\rightarrow \infty}\frac{\log\|\sigma_{s}\|_
{X\rightarrow X}}{\log s}.\]
Note that $0\leq\alpha_X\leq\beta_X\leq 1$.
In \cite[Theorem 1.5]{As} the embedding $\mathcal{M}(X)\subseteq L_{\alpha_X^{-1}}$ is shown.
The Boyd index is discussed in \cite{LT} for sequence spaces and in \cite{Ar} for s.n.\ ideals.
The Boyd index of an s.n.\ ideal $\mathfrak{S}_{\Phi}$ is defined by
\[ p=\lim_{n\rightarrow\infty}\frac{\log n}{\log\| P_{n}\|_{\Phi}} \]
when the limit exists (the limit is in $[1,\infty]$).
Theorem \ref{p-Phi-inequality} means that if $\Phi$ satisfies $(\ast)$, 
then $\mathcal{M}(\mathfrak{S}_{\Phi})\subset\mathfrak{S}_{p}$.
\end{rem}

In the rest of this paper we examine the condition ($\ast$) for a few classes of s.n.\ functions.

\begin{thm}\label{binormalize-ast}
Let $\pi$ be a binormalizing sequence and let $S_{n}$ be the partial sum defined by
$S_{n}=\sum_{j=1}^{n}\pi_{j}$.
Then $\Phi_{\pi}$ satisfies $(\ast)$ if and only if there is a constant 
$c>0$ such that for any $m,n\in\mathbb{N}$, the inequality 
\[\frac{S_{m n}}{S_{m}S_{n}}\leq c\]
holds.
\end{thm}

\begin{proof}
Let $x\in F(K)_+$ and we write its spectral decomposition by
\[ x=\sum_{j=1}^{n}s_{j}(x)p_{j}. \]
We can represent $\Phi_{\pi}(x)$ in the form
\begin{eqnarray*}
\Phi_{\pi}(x) &=& \sum_{j=1}^{n}\pi_{j}s_{j}(x)\\
    &=& (s_{1}(x)-s_{2}(x))S_{1}+\ldots +(s_{n-1}(x)-s_{n}(x))S_{n-1}+s_{n}(x)S_{n},
\end{eqnarray*}
so that if we let $e_{j}$ be the partial sum of $p_{i}$'s given by
$e_{j}=\sum_{i=1}^{j}p_{i}$, then 
\[x=(s_{1}(x)-s_{2}(x))e_{1}+\ldots +(s_{n-1}(x)-s_{n}(x))e_{n-1}+s_{n}(x)e_{n}.\]
Hence for any $a\in F(K)$,
\begin{eqnarray*}
\| x\otimes a\|_{\pi} &\leq& \bigg(\sum_{j=1}^{n}\displaystyle(s_{j}(x)-s_{j+1}(x))S_{j}\bigg)
\Max_{j}\left\{\frac{\| e_{j}\otimes a\|_{\pi}}{S_{j}}\right\}\\
    &\leq& \| x\|_{\pi}\Max_{j}\left\{\frac{\| e_{j}\otimes a\|_{\pi}}{\| e_{j}
\|_{\pi}}\right\}.
\end{eqnarray*}
Similar argument for $a$ yields
\[ \sup_{x,a}\frac{\| x\otimes a\|_{\pi}}{\| x\|_{\pi}\| a\|_{\pi}}=\sup_{p,a}
\frac{\| p\otimes a\|_{\pi}}{\| p\|_{\pi}\| a\|_{\pi}}=\sup_{p,q}
\frac{\| p\otimes q\|_{\pi}}{\| p\|_{\pi}\| q\|_{\pi}}, \]
where $p$ and $q$ run over all finite rank projections. If $p$ is a rank $n$ 
projection, then $\| p\|_{\pi}=S_{n}$ and therefore $(\ast)$ holds if and only 
if $S_{m n}/S_{m}S_{n}\leq c$.
\end{proof}

\begin{rem}
The condition 
\[ \sup_{m,n}\frac{S_{m n}}{S_{m}S_{n}}<\infty \]
appears in \cite[Theorem 6]{ACL}, as a necessarily and sufficient condition for 
the existence of exactly two nonequivalent symmetric basic sequences 
in Lorentz sequence spaces.
\end{rem}

Next we look out the Lorentz ideals $S_{p,q}$ for $1\leq q\leq p<\infty$.
When $q=1$, the Lorentz ideal $S_{p,1}$ is equal to the ideal $\mathfrak{S}_{\Phi_{\pi}}$ 
with $\pi_{j}=j^{1/p-1}$ and thus satisfies $(\ast)$ with $c_{1}=1$ 
from Theorem \ref{binormalize-ast}.

\begin{prop}
When $1\leq q\leq p<\infty$ the Lorentz ideal $S_{p,q}$ 
satisfies $(\ast)$.
\end{prop}

\begin{proof}
Let $x,y\in S_{p,q}$ be positive elements. Note that the spectrum 
of $x\otimes y$ is equal to $\{s_{i}(x)s_{j}(y)\}_{i,j=1}^{\infty}$ as a set 
considering multiplicity and each eigenspace is finite-dimensional. 
We give the  product set $\mathbb{N}\times\mathbb{N}$ an order $\prec$ by
\[(m_{1},n_{1})\prec (m_{2},n_{2}) \Longleftrightarrow \begin{cases}
    m_{1}+n_{1}<m_{2}+n_{2} \\
    \text{or} \\
    m_{1}+n_{1}=m_{2}+n_{2} \text{ and } m_{1}>m_{2}.
    \end{cases}\]
For each eigenvalue $\alpha$ of $x\otimes y$ with index $k$, let $I_{\alpha}$ 
be the finite sequence $\{(m_{1},n_{1}),\ldots,(m_k,n_k)\}$ 
in $\mathbb{N}\times\mathbb{N}$ such that 
$s_{m_{i}}(x)s_{n_{i}}(y)=\alpha$ and $(m_{i},n_{i})\prec(m_{i+1},n_{i+1})$. 
If $s_{j+1}(x\otimes y)=\cdots=s_{j+k}(x\otimes y)=\alpha$, 
for all $i=1,\ldots,k$ we have 
\[s_{j+i}(x\otimes y)=s_{m_{i}}(x)s_{n_{i}}(y)\]
and $j+i\geq m_{i}n_{i}$.
Hence 
\begin{eqnarray*}
\| x\otimes y\|_{p,q} &=& \left(\sum_{j=1}^{\infty}
    \frac{s_{j}(x\otimes y)^q}{j^{1-q/p}}\right)^{1/q}\\
 &\leq& \left(\sum_{i,j=1}^{\infty}
    \frac{s_{i}(x)^{q}s_{j}(y)^{q}}{(i j)^{1-q/p}}\right)^{1/q}\\
 &=& \left(\sum_{j=1}^{\infty}\frac{s_{j}(x)^{q}}{j^{1-q/p}}\right)^{1/q}
 \left(\sum_{j=1}^{\infty}\frac{s_{j}(y)^{q}}{j^{1-q/p}}\right)^{1/q}
 =\| x\|_{p,q}\| y\|_{p,q}.
\end{eqnarray*}
\end{proof}

\begin{rem}
In \cite[p.253]{AMS} it is shown that for the Lorentz function space 
$L_{p,q}~(1<p<\infty,~1\leq q\leq\infty)$, 
we have $\mathcal{M}(L_{p,q})=L_{p,\min(p,q)}$.
\end{rem}

\section*{Acknowledgement}
The author would like to thank M.\ Izumi for suggesting this problem.

\end{document}